\documentclass{amsproc}
\usepackage{color}
\newtheorem{theorem}{Theorem}[section]
\newtheorem{lemma}[theorem]{Lemma}

\newtheorem{corollary}[theorem]{Corollary}

\theoremstyle{definition}
\newtheorem{definition}[theorem]{Definition}
\newtheorem{example}[theorem]{Example}

\newtheorem{question}[theorem]{Question}

\theoremstyle{remark}
\newtheorem{remark}[theorem]{Remark}

\numberwithin{equation}{section}

\newcommand{\abs}[1]{\lvert#1\rvert}

\newcommand{\bd}{{\sc Proof}.\ \ }

\begin{document}

\title
{Horocycle flows without minimal sets}

\author{Shigenori Matsumoto}
\address{Department of Mathematics, College of
Science and Technology, Nihon University, 1-8-14 Kanda, Surugadai,
Chiyoda-ku, Tokyo, 101-8308 Japan}
\email{matsumo@math.cst.nihon-u.ac.jp}

\thanks{The author is partially supported by Grant-in-Aid for
Scientific Research (C) No.\ 25400096.}

\date{\today}

\newcommand{\AAA}{{\mathbb A}}
\newcommand{\BBB}{{\mathbb B}}
\newcommand{\LL}{{\mathcal L}}
\newcommand{\MCG}{{\rm MCG}}
\newcommand{\PSL}{{\rm PSL}}
\newcommand{\R}{{\mathbb R}}
\newcommand{\Z}{{\mathbb Z}}
\newcommand{\XX}{{\mathcal X}}
\newcommand{\per}{{\rm per}}
\newcommand{\N}{{\mathbb N}}

\newcommand{\PP}{{\mathcal P}}
\newcommand{\GG}{{\mathbb G}}
\newcommand{\FF}{{\mathcal F}}
\newcommand{\EE}{{\mathbb E}}
\newcommand{\BB}{{\mathbb B}}
\newcommand{\CC}{{\mathcal C}}
\newcommand{\HH}{{\mathbb H}}
\newcommand{\UU}{{\mathcal U}}
\newcommand{\oboundary}{{\mathbb S}^1_\infty}
\newcommand{\Q}{{\mathbb Q}}
\newcommand{\DD}{{\mathcal D}}
\newcommand{\rot}{{\rm rot}}
\newcommand{\Cl}{{\rm Cl}}
\newcommand{\Index}{{\rm Index}}
\newcommand{\Int}{{\rm Int}}
\newcommand{\Fr}{{\rm Fr}}

\def\bm{ \left[ \begin{array}{ccc} }
\def\ema{\end{array} \right] }
\def\bmt{ \left[ \begin{array}{cc} }

\date{\today }

\maketitle

\begin{abstract}
 We show that the horocycle flows of open tight hyperbolic surfaces
do not admit minimal sets.
\end{abstract}

\section{Introduction}

Let $\{\phi^t\}$ be a flow of a metric space $X$.  A subset 
of $X$ is called a {\em minimal set} of $\{\phi^t\}$ if it is
closed and invariant by $\phi^t$, and is minimal among them
with respect to the inclusion.
 If $X$ is compact, then any flow on $X$
admits a minimal set. But if $X$ is not compact, this is
not always the case. The first example of a flow without minimal set is
constructed on an open surface by T. Inaba \cite{I}. Later various
examples are piled up by many authors including \cite{BM}. See also
\cite{MS} for examples of Anzai skew products on an open annulus.

M. Kulikov \cite{K} constructed an example of the horocycle flow
of an open hyperbolic surface with this property.
This is interesting since horocycle flows have long been
studied by various mathematicians; function analysists, topologists,
dynamical people and ergodic theoretists. Moreover an example in
\cite{K} is the first one constructed algebraically on an homogeneous
space. However the example is constructed in
a specific and elaborate way. The purpose of this paper is to show
that more general Fuchsian groups also satify this property. 

\begin{definition} \label{def}
 A Fuchsian group $\Gamma$ is called {\em tight} if it satisfies the
 following conditions.

(1) $\Gamma$ is purely hyperbolic.

(2)
$\Sigma=\Gamma\setminus\HH^2$ is noncompact and admits an
increasing and exausting  sequence of compact subsurfaces $\{\Sigma_n\}_{n\in\N}$ with geodesic
 boundaries such that there is a bound $C$ on the length of
 components of $\partial \Sigma_n$.
\end{definition}

Tight Fuchsian groups are infinitely generated and of the first kind.
The main result of this notes is  the following.

\begin{theorem}\label{main}
If $\Gamma$ is a tight Fuchsian group, the horocycle flow $\{h^s\}_{s\in\R}$
on $\Gamma\setminus PSL(2,\R)$ admits no minimal sets.
\end{theorem} 

\begin{corollary}\label{cor}
 Almost all orbits of the horocycle flow of a tight Fuchsian group are dense.
\end{corollary}

\begin{remark} \label{rem}
 For some tight Fuchsian groups, the horocycle flow is ergodic, while 
for others it is not.
\end{remark}

In Section 2, we explain conventions used in this paper. In Section 3, we
prepare fundamental facts about the horocyclic limit points.
Section 4 is devoted to the proof of Theorem \ref{main}. Finally in Section
5, we raise examples of tight Fuchsian groups and discuss 
Corollary \ref{cor} and Remark \ref{rem}.

\section{Conventions}

The right coset space $PSL(2,\R)/PSO(2)$ is identified with the upper
half plane $\HH$ by sending a matrix $\pm\bmt a&b\\c&d\ema$ to a point
$\displaystyle \frac{ai+b}{ci+d}$. 
The group $PSL(2,\R)$ acts on $\HH$ as linear fractional
transformations,
and is identified
with the unit tangent space $T^1\HH$ by sending $M\in PSL(2,\R)$
to $M_*(i,\vec e\,)$, where $(i,\vec e\,)$ is the upward unit tangent
vector at $i$. The canonical projection is denoted by
$$\pi_1:PSL(2,\R)=T^1\HH\to\HH.$$

The geodesic flow $\{\tilde g^t\}$ (resp.\ the horocycle flow $\{\tilde h^s\}$)
 on $PSL(2,\R)$ is given by the right multiplication
of the matrices $\bmt e^{t/2}&0\\0&e^{-t/2}\ema$ (resp.\ $\bmt 1 &s\\0&1\ema$).
The quotient space $PSL(2,\R)/\langle \tilde h^s\rangle$ is identified with
the annulus $\AAA=(\R^2\setminus\{0\})/\langle\pm1\rangle$, by sending
$M\in PSL(2,\R)$ to a point $M(1,0)^t\in\AAA$.
The canonical projection is denoted by 
$$\pi_2:T^1\HH=PSL(2,\R)\to\AAA.$$
The geodesic flow $\tilde g^t$ induces a flow on $\AAA$, which is just
the scalar multiplication by $e^{t/2}$. The further quotient space
$\AAA/\R_+$ is equal to
the circle at infinity $\partial_\infty\HH$, both being defined as
the right coset space of $PSL(2,\R)$ by the subgroup of the upper triangular
 matrices. 

For any $\xi\in
\partial_\infty\HH$, the preimage of $\xi$ by the canonical projection
$\AAA\to \partial_\infty\HH$ is denoted by $\AAA(\xi)$. 
It is a ray of $\AAA$.
For any point $p\in \AAA(\xi)$, $H(p)=\pi_1(\pi_2^{-1}(p))$
is a horocycle in $\HH$ tangent to $\partial_\infty\HH$ at $\xi$.
The open horodisk encircled by $H(p)$ is denoted by $D(p)$. 
The signed distance from $i\in\HH^2$
to the horocycle $H(p)$ (positive if $i$ is outside $D(p)$ and negative
if inside) is  $2\log\abs{p}$,
where $\abs{p}$ denotes the Euclidian norm of $\AAA$.
Thus $\abs{p}<1$ if and only if $i\in D(p)$.

Given a Fuchsian group $\Gamma$, 
the flows $\{\tilde g^t\}$ and $\{\tilde h^s\}$
induce flows on $\Gamma\setminus PSL(2,\R)$,  denoted
by $\{g^t\}$ and $\{h^s\}$. The right $\R$-action $\{h^s\}$ on 
$\Gamma\setminus PSL(2,\R)$
is Morita equivalent to the left $\Gamma$-action on $\AAA$. Thus a dense 
$\{h^s\}$-orbit in $\Gamma\setminus PSL(2,\R)$ corresponds to a 
dense $\Gamma$-orbit in $\AAA$. Likewise a minimal set of the flow $\{h^s\}$
in $\Gamma\setminus PSL(2,\R)$ corresponds to a minimal set for the $\Gamma$-action on
$\AAA$. 


For any $\tilde v\in T^1\HH$,
$t\mapsto\pi_1\tilde g^t(\tilde v)$ is the unit speed
geodesic in $\HH$ with 
innitial vector $\tilde v$. Its positive
endpoint in $\partial_\infty\HH$ is denoted by
$\tilde v(\infty)$. 
If $\Gamma$ is purely hyperbolic, the quotient space 
$\Sigma=\Gamma\setminus\HH$ is a hyperbolic surface,
and its unit tangent bundle $T^1\Sigma$ is identified with
$\Gamma\setminus PSL(2,\R)$.
The canonical projection is denoted by 
$$\pi:T^1\Sigma\to \Sigma.$$
 For any
$v\in T^1\Sigma$, 
$$v[0,\infty)=\{\pi g^t(v)\mid 0\leq t<\infty\}$$ 
is the geodesic ray in $\Sigma$
with innitial vector $v$. 

\section{Horocyclic limit points}
In this section we assume that $\Gamma$ is a purely hyperbolic Fuchsian
group of the first kind. As before, we denote $\Sigma=\Gamma\setminus\HH$.
Many of the contents in this section are taken from \cite{S}.

\begin{definition}
A geodesic ray $v[0,\infty)$ , $v\in T^1\Sigma$, is called a {\em
quasi-minimizer} if there is $k>0$ such that
$d(\pi g^t(v),\pi(v))\geq t-k$ for any $t\geq0$.
\end{definition}

See Figure 1.
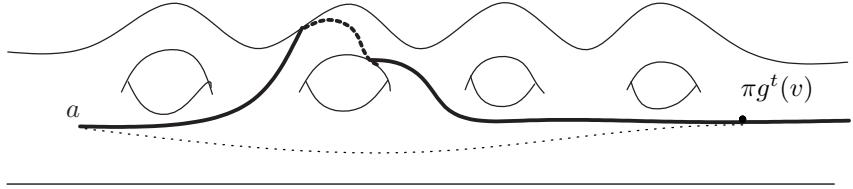
\begin{figure}[h]
{\unitlength 0.1in%
\begin{picture}( 44.0700,  9.5000)(  7.0000,-23.7000)%
%
\special{pn 8}%
\special{pa 704 1678}%
\special{pa 736 1681}%
\special{pa 769 1683}%
\special{pa 801 1685}%
\special{pa 834 1687}%
\special{pa 866 1688}%
\special{pa 898 1688}%
\special{pa 930 1687}%
\special{pa 962 1685}%
\special{pa 994 1682}%
\special{pa 1026 1678}%
\special{pa 1058 1672}%
\special{pa 1089 1664}%
\special{pa 1120 1655}%
\special{pa 1151 1643}%
\special{pa 1182 1629}%
\special{pa 1212 1613}%
\special{pa 1242 1596}%
\special{pa 1272 1577}%
\special{pa 1302 1557}%
\special{pa 1389 1497}%
\special{pa 1418 1479}%
\special{pa 1446 1463}%
\special{pa 1475 1449}%
\special{pa 1503 1437}%
\special{pa 1531 1429}%
\special{pa 1559 1424}%
\special{pa 1587 1423}%
\special{pa 1614 1428}%
\special{pa 1642 1436}%
\special{pa 1669 1449}%
\special{pa 1696 1466}%
\special{pa 1723 1485}%
\special{pa 1751 1506}%
\special{pa 1805 1550}%
\special{pa 1832 1573}%
\special{pa 1860 1594}%
\special{pa 1888 1613}%
\special{pa 1915 1630}%
\special{pa 1943 1644}%
\special{pa 1972 1654}%
\special{pa 2000 1660}%
\special{pa 2029 1660}%
\special{pa 2059 1655}%
\special{pa 2088 1645}%
\special{pa 2118 1632}%
\special{pa 2147 1616}%
\special{pa 2177 1598}%
\special{pa 2207 1578}%
\special{pa 2237 1557}%
\special{pa 2267 1535}%
\special{pa 2297 1514}%
\special{pa 2327 1494}%
\special{pa 2356 1475}%
\special{pa 2386 1459}%
\special{pa 2415 1445}%
\special{pa 2443 1435}%
\special{pa 2472 1429}%
\special{pa 2500 1428}%
\special{pa 2527 1433}%
\special{pa 2554 1443}%
\special{pa 2581 1459}%
\special{pa 2607 1478}%
\special{pa 2633 1500}%
\special{pa 2659 1524}%
\special{pa 2685 1549}%
\special{pa 2710 1574}%
\special{pa 2736 1597}%
\special{pa 2762 1619}%
\special{pa 2788 1637}%
\special{pa 2814 1650}%
\special{pa 2840 1659}%
\special{pa 2867 1661}%
\special{pa 2894 1658}%
\special{pa 2922 1649}%
\special{pa 2950 1636}%
\special{pa 2978 1619}%
\special{pa 3006 1600}%
\special{pa 3035 1578}%
\special{pa 3063 1556}%
\special{pa 3092 1532}%
\special{pa 3121 1509}%
\special{pa 3150 1487}%
\special{pa 3179 1467}%
\special{pa 3209 1450}%
\special{pa 3238 1435}%
\special{pa 3267 1425}%
\special{pa 3296 1420}%
\special{pa 3324 1421}%
\special{pa 3353 1427}%
\special{pa 3382 1439}%
\special{pa 3410 1454}%
\special{pa 3438 1474}%
\special{pa 3466 1495}%
\special{pa 3550 1567}%
\special{pa 3577 1590}%
\special{pa 3604 1612}%
\special{pa 3631 1631}%
\special{pa 3658 1647}%
\special{pa 3685 1659}%
\special{pa 3712 1665}%
\special{pa 3738 1666}%
\special{pa 3764 1660}%
\special{pa 3790 1649}%
\special{pa 3816 1633}%
\special{pa 3842 1613}%
\special{pa 3868 1591}%
\special{pa 3920 1541}%
\special{pa 3945 1517}%
\special{pa 3971 1493}%
\special{pa 3998 1471}%
\special{pa 4024 1452}%
\special{pa 4051 1438}%
\special{pa 4077 1428}%
\special{pa 4105 1423}%
\special{pa 4132 1423}%
\special{pa 4160 1427}%
\special{pa 4187 1435}%
\special{pa 4216 1446}%
\special{pa 4244 1461}%
\special{pa 4273 1477}%
\special{pa 4301 1496}%
\special{pa 4331 1516}%
\special{pa 4360 1537}%
\special{pa 4389 1559}%
\special{pa 4419 1581}%
\special{pa 4479 1623}%
\special{pa 4510 1642}%
\special{pa 4540 1660}%
\special{pa 4571 1676}%
\special{pa 4602 1689}%
\special{pa 4633 1701}%
\special{pa 4664 1711}%
\special{pa 4696 1719}%
\special{pa 4727 1726}%
\special{pa 4759 1732}%
\special{pa 4790 1736}%
\special{pa 4854 1740}%
\special{pa 4886 1741}%
\special{pa 4918 1741}%
\special{pa 4950 1740}%
\special{pa 4983 1739}%
\special{pa 5015 1738}%
\special{pa 5047 1736}%
\special{pa 5079 1733}%
\special{pa 5100 1732}%
\special{fp}%
%
\special{pn 8}%
\special{pa 1299 1888}%
\special{pa 1313 1858}%
\special{pa 1328 1828}%
\special{pa 1343 1799}%
\special{pa 1361 1772}%
\special{pa 1380 1747}%
\special{pa 1403 1724}%
\special{pa 1428 1705}%
\special{pa 1457 1689}%
\special{pa 1489 1677}%
\special{pa 1523 1669}%
\special{pa 1557 1666}%
\special{pa 1590 1667}%
\special{pa 1622 1673}%
\special{pa 1650 1685}%
\special{pa 1675 1701}%
\special{pa 1696 1722}%
\special{pa 1714 1747}%
\special{pa 1730 1776}%
\special{pa 1743 1807}%
\special{pa 1756 1840}%
\special{pa 1767 1874}%
\special{pa 1775 1900}%
\special{fp}%
%
\special{pn 8}%
\special{pa 1334 1834}%
\special{pa 1348 1864}%
\special{pa 1364 1894}%
\special{pa 1380 1921}%
\special{pa 1399 1947}%
\special{pa 1421 1969}%
\special{pa 1447 1987}%
\special{pa 1475 1999}%
\special{pa 1507 2005}%
\special{pa 1541 2004}%
\special{pa 1576 1994}%
\special{pa 1611 1974}%
\special{pa 1647 1944}%
\special{pa 1681 1908}%
\special{pa 1712 1873}%
\special{pa 1738 1846}%
\special{pa 1757 1833}%
\special{pa 1767 1839}%
\special{pa 1768 1864}%
\special{fp}%
%
\special{pn 8}%
\special{pa 2230 1888}%
\special{pa 2246 1859}%
\special{pa 2263 1830}%
\special{pa 2281 1802}%
\special{pa 2300 1776}%
\special{pa 2322 1753}%
\special{pa 2345 1733}%
\special{pa 2372 1717}%
\special{pa 2401 1705}%
\special{pa 2433 1697}%
\special{pa 2466 1693}%
\special{pa 2500 1693}%
\special{pa 2533 1699}%
\special{pa 2566 1708}%
\special{pa 2598 1723}%
\special{pa 2627 1742}%
\special{pa 2654 1765}%
\special{pa 2676 1791}%
\special{pa 2692 1819}%
\special{pa 2703 1850}%
\special{pa 2706 1881}%
\special{pa 2706 1882}%
\special{fp}%
%
\special{pn 8}%
\special{pa 2265 1834}%
\special{pa 2280 1864}%
\special{pa 2296 1893}%
\special{pa 2314 1920}%
\special{pa 2335 1943}%
\special{pa 2360 1962}%
\special{pa 2389 1976}%
\special{pa 2420 1985}%
\special{pa 2454 1988}%
\special{pa 2488 1987}%
\special{pa 2521 1980}%
\special{pa 2552 1969}%
\special{pa 2579 1954}%
\special{pa 2604 1934}%
\special{pa 2627 1912}%
\special{pa 2648 1887}%
\special{pa 2669 1861}%
\special{pa 2688 1834}%
\special{pa 2692 1828}%
\special{fp}%
%
\special{pn 8}%
\special{pa 3098 1876}%
\special{pa 3110 1845}%
\special{pa 3124 1814}%
\special{pa 3139 1786}%
\special{pa 3158 1760}%
\special{pa 3181 1738}%
\special{pa 3208 1721}%
\special{pa 3240 1709}%
\special{pa 3274 1703}%
\special{pa 3309 1702}%
\special{pa 3342 1708}%
\special{pa 3373 1719}%
\special{pa 3399 1736}%
\special{pa 3420 1759}%
\special{pa 3438 1786}%
\special{pa 3454 1814}%
\special{pa 3471 1843}%
\special{pa 3490 1869}%
\special{pa 3511 1894}%
\special{fp}%
%
\special{pn 8}%
\special{pa 3126 1804}%
\special{pa 3136 1837}%
\special{pa 3148 1868}%
\special{pa 3162 1897}%
\special{pa 3181 1922}%
\special{pa 3205 1942}%
\special{pa 3234 1956}%
\special{pa 3267 1963}%
\special{pa 3301 1964}%
\special{pa 3334 1959}%
\special{pa 3364 1948}%
\special{pa 3390 1930}%
\special{pa 3414 1908}%
\special{pa 3435 1883}%
\special{pa 3455 1856}%
\special{pa 3462 1846}%
\special{fp}%
%
\special{pn 8}%
\special{pa 3945 1876}%
\special{pa 3961 1846}%
\special{pa 3978 1817}%
\special{pa 3997 1791}%
\special{pa 4018 1768}%
\special{pa 4043 1750}%
\special{pa 4073 1738}%
\special{pa 4105 1731}%
\special{pa 4139 1731}%
\special{pa 4173 1736}%
\special{pa 4205 1746}%
\special{pa 4234 1761}%
\special{pa 4258 1781}%
\special{pa 4279 1805}%
\special{pa 4298 1832}%
\special{pa 4315 1860}%
\special{pa 4330 1888}%
\special{fp}%
%
\special{pn 8}%
\special{pa 3966 1834}%
\special{pa 3977 1868}%
\special{pa 3990 1900}%
\special{pa 4005 1928}%
\special{pa 4024 1951}%
\special{pa 4049 1967}%
\special{pa 4079 1975}%
\special{pa 4112 1976}%
\special{pa 4147 1972}%
\special{pa 4180 1961}%
\special{pa 4211 1947}%
\special{pa 4237 1928}%
\special{pa 4261 1907}%
\special{pa 4283 1883}%
\special{pa 4303 1859}%
\special{pa 4309 1852}%
\special{fp}%
%
\special{pn 20}%
\special{pa 1082 2068}%
\special{pa 1114 2069}%
\special{pa 1145 2069}%
\special{pa 1177 2070}%
\special{pa 1208 2070}%
\special{pa 1240 2071}%
\special{pa 1303 2071}%
\special{pa 1335 2070}%
\special{pa 1367 2070}%
\special{pa 1399 2069}%
\special{pa 1432 2068}%
\special{pa 1464 2067}%
\special{pa 1530 2063}%
\special{pa 1596 2057}%
\special{pa 1629 2053}%
\special{pa 1695 2043}%
\special{pa 1728 2037}%
\special{pa 1760 2030}%
\special{pa 1792 2022}%
\special{pa 1823 2012}%
\special{pa 1853 2002}%
\special{pa 1882 1990}%
\special{pa 1911 1977}%
\special{pa 1938 1962}%
\special{pa 1964 1946}%
\special{pa 1989 1928}%
\special{pa 2013 1908}%
\special{pa 2035 1887}%
\special{pa 2056 1865}%
\special{pa 2076 1841}%
\special{pa 2096 1816}%
\special{pa 2114 1790}%
\special{pa 2132 1763}%
\special{pa 2149 1735}%
\special{pa 2181 1677}%
\special{pa 2196 1647}%
\special{pa 2212 1617}%
\special{pa 2242 1557}%
\special{pa 2244 1552}%
\special{fp}%
%
\special{pn 20}%
\special{pn 20}%
\special{pa 2258 1552}%
\special{pa 2275 1543}%
\special{fp}%
\special{pa 2297 1533}%
\special{pa 2315 1525}%
\special{fp}%
\special{pa 2338 1517}%
\special{pa 2357 1513}%
\special{fp}%
\special{pa 2381 1510}%
\special{pa 2401 1511}%
\special{fp}%
\special{pa 2425 1515}%
\special{pa 2444 1520}%
\special{fp}%
\special{pa 2466 1530}%
\special{pa 2483 1539}%
\special{fp}%
\special{pa 2503 1553}%
\special{pa 2516 1568}%
\special{fp}%
\special{pa 2530 1587}%
\special{pa 2540 1604}%
\special{fp}%
\special{pa 2549 1627}%
\special{pa 2556 1645}%
\special{fp}%
\special{pa 2565 1668}%
\special{pa 2574 1685}%
\special{fp}%
\special{pa 2587 1706}%
\special{pa 2599 1721}%
\special{fp}%
\special{pa 2619 1736}%
\special{pa 2636 1744}%
\special{fp}%
%
\special{pn 20}%
\special{pa 2608 1720}%
\special{pa 2640 1721}%
\special{pa 2672 1723}%
\special{pa 2704 1727}%
\special{pa 2735 1734}%
\special{pa 2766 1744}%
\special{pa 2795 1758}%
\special{pa 2823 1774}%
\special{pa 2849 1793}%
\special{pa 2873 1815}%
\special{pa 2895 1840}%
\special{pa 2915 1866}%
\special{pa 2934 1892}%
\special{pa 2953 1919}%
\special{pa 2973 1943}%
\special{pa 2993 1966}%
\special{pa 3016 1985}%
\special{pa 3042 2001}%
\special{pa 3069 2014}%
\special{pa 3098 2023}%
\special{pa 3128 2031}%
\special{pa 3160 2036}%
\special{pa 3193 2039}%
\special{pa 3227 2041}%
\special{pa 3262 2042}%
\special{pa 3402 2038}%
\special{pa 3436 2036}%
\special{pa 3471 2036}%
\special{pa 3505 2035}%
\special{pa 3538 2034}%
\special{pa 3734 2034}%
\special{pa 3765 2035}%
\special{pa 3797 2035}%
\special{pa 3829 2036}%
\special{pa 3860 2036}%
\special{pa 3922 2038}%
\special{pa 3953 2038}%
\special{pa 4046 2041}%
\special{pa 4077 2041}%
\special{pa 4139 2043}%
\special{pa 4170 2043}%
\special{pa 4201 2044}%
\special{pa 4232 2044}%
\special{pa 4264 2045}%
\special{pa 4326 2045}%
\special{pa 4358 2046}%
\special{pa 4548 2046}%
\special{pa 4580 2045}%
\special{pa 4677 2045}%
\special{pa 4709 2044}%
\special{pa 4741 2044}%
\special{pa 4773 2043}%
\special{pa 4838 2043}%
\special{pa 4870 2042}%
\special{pa 4902 2042}%
\special{pa 4935 2041}%
\special{pa 4967 2040}%
\special{pa 5000 2040}%
\special{pa 5032 2039}%
\special{pa 5064 2039}%
\special{pa 5097 2038}%
\special{pa 5107 2038}%
\special{fp}%
\put(10.1000,-20.2000){\makebox(0,0)[lb]{$a$}}%
%
\special{pn 4}%
\special{sh 1}%
\special{ar 4550 2030 16 16 0  6.28318530717959E+0000}%
\special{sh 1}%
\special{ar 4550 2030 16 16 0  6.28318530717959E+0000}%
\special{sh 1}%
\special{ar 4550 2030 16 16 0  6.28318530717959E+0000}%
\put(45.4000,-19.4000){\makebox(0,0)[lb]{$\pi g^t(v)$}}%
%
\special{pn 8}%
\special{pa 700 2370}%
\special{pa 5030 2370}%
\special{fp}%
%
\special{pn 8}%
\special{pn 8}%
\special{pa 1110 2080}%
\special{pa 1118 2081}%
\special{fp}%
\special{pa 1155 2086}%
\special{pa 1163 2087}%
\special{fp}%
\special{pa 1200 2091}%
\special{pa 1208 2092}%
\special{fp}%
\special{pa 1245 2097}%
\special{pa 1253 2098}%
\special{fp}%
\special{pa 1290 2102}%
\special{pa 1298 2103}%
\special{fp}%
\special{pa 1335 2108}%
\special{pa 1343 2109}%
\special{fp}%
\special{pa 1380 2113}%
\special{pa 1388 2114}%
\special{fp}%
\special{pa 1425 2118}%
\special{pa 1433 2119}%
\special{fp}%
\special{pa 1470 2124}%
\special{pa 1478 2125}%
\special{fp}%
\special{pa 1515 2130}%
\special{pa 1523 2131}%
\special{fp}%
\special{pa 1560 2134}%
\special{pa 1568 2135}%
\special{fp}%
\special{pa 1605 2140}%
\special{pa 1613 2140}%
\special{fp}%
\special{pa 1650 2145}%
\special{pa 1658 2145}%
\special{fp}%
\special{pa 1695 2149}%
\special{pa 1703 2150}%
\special{fp}%
\special{pa 1740 2154}%
\special{pa 1748 2155}%
\special{fp}%
\special{pa 1785 2159}%
\special{pa 1793 2160}%
\special{fp}%
\special{pa 1830 2164}%
\special{pa 1838 2164}%
\special{fp}%
\special{pa 1876 2168}%
\special{pa 1884 2169}%
\special{fp}%
\special{pa 1921 2172}%
\special{pa 1929 2172}%
\special{fp}%
\special{pa 1966 2175}%
\special{pa 1974 2176}%
\special{fp}%
\special{pa 2011 2179}%
\special{pa 2019 2180}%
\special{fp}%
\special{pa 2056 2183}%
\special{pa 2064 2184}%
\special{fp}%
\special{pa 2102 2186}%
\special{pa 2110 2187}%
\special{fp}%
\special{pa 2147 2189}%
\special{pa 2155 2190}%
\special{fp}%
\special{pa 2192 2193}%
\special{pa 2200 2193}%
\special{fp}%
\special{pa 2237 2195}%
\special{pa 2245 2196}%
\special{fp}%
\special{pa 2282 2197}%
\special{pa 2291 2198}%
\special{fp}%
\special{pa 2328 2199}%
\special{pa 2336 2200}%
\special{fp}%
\special{pa 2373 2202}%
\special{pa 2381 2202}%
\special{fp}%
\special{pa 2418 2203}%
\special{pa 2426 2203}%
\special{fp}%
\special{pa 2464 2204}%
\special{pa 2472 2205}%
\special{fp}%
\special{pa 2509 2206}%
\special{pa 2517 2206}%
\special{fp}%
\special{pa 2554 2206}%
\special{pa 2562 2206}%
\special{fp}%
\special{pa 2600 2207}%
\special{pa 2608 2207}%
\special{fp}%
\special{pa 2645 2207}%
\special{pa 2653 2207}%
\special{fp}%
\special{pa 2690 2207}%
\special{pa 2698 2206}%
\special{fp}%
\special{pa 2736 2206}%
\special{pa 2744 2206}%
\special{fp}%
\special{pa 2781 2205}%
\special{pa 2789 2204}%
\special{fp}%
\special{pa 2826 2203}%
\special{pa 2834 2203}%
\special{fp}%
\special{pa 2872 2202}%
\special{pa 2880 2202}%
\special{fp}%
\special{pa 2917 2200}%
\special{pa 2925 2199}%
\special{fp}%
\special{pa 2962 2197}%
\special{pa 2970 2196}%
\special{fp}%
\special{pa 3007 2194}%
\special{pa 3015 2194}%
\special{fp}%
\special{pa 3053 2191}%
\special{pa 3061 2191}%
\special{fp}%
\special{pa 3098 2187}%
\special{pa 3106 2187}%
\special{fp}%
\special{pa 3143 2184}%
\special{pa 3151 2183}%
\special{fp}%
\special{pa 3188 2180}%
\special{pa 3196 2179}%
\special{fp}%
\special{pa 3233 2175}%
\special{pa 3241 2175}%
\special{fp}%
\special{pa 3278 2170}%
\special{pa 3286 2169}%
\special{fp}%
\special{pa 3323 2166}%
\special{pa 3331 2165}%
\special{fp}%
\special{pa 3368 2161}%
\special{pa 3376 2160}%
\special{fp}%
\special{pa 3414 2156}%
\special{pa 3422 2155}%
\special{fp}%
\special{pa 3458 2151}%
\special{pa 3466 2150}%
\special{fp}%
\special{pa 3504 2146}%
\special{pa 3512 2145}%
\special{fp}%
\special{pa 3549 2140}%
\special{pa 3557 2139}%
\special{fp}%
\special{pa 3594 2135}%
\special{pa 3602 2134}%
\special{fp}%
\special{pa 3639 2130}%
\special{pa 3647 2129}%
\special{fp}%
\special{pa 3684 2124}%
\special{pa 3692 2123}%
\special{fp}%
\special{pa 3729 2119}%
\special{pa 3737 2118}%
\special{fp}%
\special{pa 3774 2114}%
\special{pa 3782 2113}%
\special{fp}%
\special{pa 3819 2108}%
\special{pa 3827 2108}%
\special{fp}%
\special{pa 3864 2104}%
\special{pa 3872 2103}%
\special{fp}%
\special{pa 3909 2099}%
\special{pa 3917 2098}%
\special{fp}%
\special{pa 3954 2095}%
\special{pa 3962 2094}%
\special{fp}%
\special{pa 3999 2090}%
\special{pa 4007 2089}%
\special{fp}%
\special{pa 4044 2085}%
\special{pa 4052 2085}%
\special{fp}%
\special{pa 4089 2082}%
\special{pa 4097 2081}%
\special{fp}%
\special{pa 4134 2078}%
\special{pa 4142 2077}%
\special{fp}%
\special{pa 4180 2074}%
\special{pa 4188 2074}%
\special{fp}%
\special{pa 4225 2071}%
\special{pa 4233 2070}%
\special{fp}%
\special{pa 4270 2068}%
\special{pa 4278 2068}%
\special{fp}%
\special{pa 4315 2066}%
\special{pa 4323 2066}%
\special{fp}%
\special{pa 4361 2064}%
\special{pa 4369 2064}%
\special{fp}%
\special{pa 4406 2062}%
\special{pa 4414 2061}%
\special{fp}%
\special{pa 4451 2061}%
\special{pa 4459 2061}%
\special{fp}%
\special{pa 4497 2060}%
\special{pa 4505 2060}%
\special{fp}%
\special{pa 4542 2060}%
\special{pa 4550 2060}%
\special{fp}%
\end{picture}}%

\caption{\small A quasi-minimizer. It starts at a point $a$ and
after one turn goes straight to the right. Any point $\pi g^t(v)$ on the
 curve satisfies $t-d(a,\pi g^t(v))\leq k$ for some $k$, where $t$ is
 the length of the curve between the two points.}
\end{figure}

\begin{definition}
 A point at infinity $\xi\in \partial_\infty\HH$ is called a {\em horocyclic
limit point} of $\Gamma$ if any horodisk at $\xi$ intersects the orbit
$\Gamma i$. Otherwise it is called {\em nonhorocyclic.}
\end{definition}

See Figure 2.
If $\xi$ is a horocyclic limit point, then any horodisk at $\xi$
intersects any orbit $\Gamma z$.
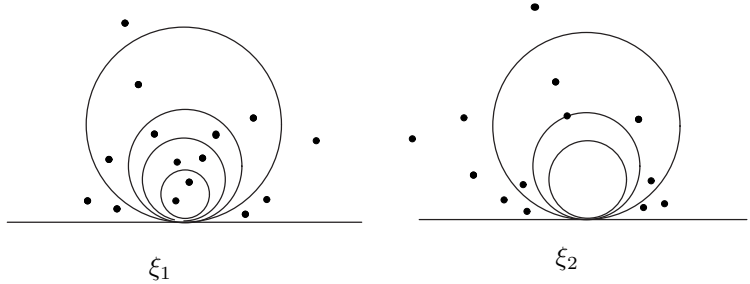
\begin{figure}[h]
{\unitlength 0.1in%
\begin{picture}( 38.7100, 13.0000)(  6.7000,-20.7000)%
%
\special{pn 8}%
\special{pa 670 1897}%
\special{pa 2525 1897}%
\special{fp}%
%
\special{pn 8}%
\special{pa 2826 1883}%
\special{pa 4541 1883}%
\special{fp}%
%
\special{pn 8}%
\special{ar 1594 1673 217 217  1.7975952  1.5707963}%
%
\special{pn 8}%
\special{ar 3708 1673 203 203  0.0000000  6.2831853}%
%
\special{pn 8}%
\special{ar 1601 1603 297 297  0.0000000  6.2831853}%
%
\special{pn 8}%
\special{ar 3701 1603 280 280  0.0000000  6.2831853}%
%
\special{pn 8}%
\special{ar 1594 1386 511 511  0.0000000  6.2831853}%
%
\special{pn 8}%
\special{ar 3701 1393 490 490  0.0000000  6.2831853}%
%
\special{pn 8}%
\special{ar 1601 1750 127 127  0.0000000  6.2831853}%
%
\special{pn 4}%
\special{sh 1}%
\special{ar 1286 854 16 16 0  6.28318530717959E+0000}%
\special{sh 1}%
\special{ar 1286 854 16 16 0  6.28318530717959E+0000}%
%
\special{pn 4}%
\special{sh 1}%
\special{ar 1356 1176 16 16 0  6.28318530717959E+0000}%
\special{sh 1}%
\special{ar 1356 1176 16 16 0  6.28318530717959E+0000}%
%
\special{pn 4}%
\special{sh 1}%
\special{ar 1440 1435 16 16 0  6.28318530717959E+0000}%
\special{sh 1}%
\special{ar 1440 1435 16 16 0  6.28318530717959E+0000}%
%
\special{pn 4}%
\special{sh 1}%
\special{ar 1692 1561 16 16 0  6.28318530717959E+0000}%
\special{sh 1}%
\special{ar 1692 1561 16 16 0  6.28318530717959E+0000}%
%
\special{pn 4}%
\special{sh 1}%
\special{ar 1622 1687 16 16 0  6.28318530717959E+0000}%
\special{sh 1}%
\special{ar 1622 1687 16 16 0  6.28318530717959E+0000}%
%
\special{pn 4}%
\special{sh 1}%
\special{ar 1552 1785 16 16 0  6.28318530717959E+0000}%
\special{sh 1}%
\special{ar 1559 1582 16 16 0  6.28318530717959E+0000}%
%
\special{pn 4}%
\special{sh 1}%
\special{ar 1202 1568 16 16 0  6.28318530717959E+0000}%
\special{sh 1}%
\special{ar 1202 1568 16 16 0  6.28318530717959E+0000}%
%
\special{pn 4}%
\special{sh 1}%
\special{ar 1958 1351 16 16 0  6.28318530717959E+0000}%
\special{sh 1}%
\special{ar 1958 1351 16 16 0  6.28318530717959E+0000}%
%
\special{pn 4}%
\special{sh 1}%
\special{ar 1762 1435 16 16 0  6.28318530717959E+0000}%
\special{sh 1}%
\special{ar 1762 1442 16 16 0  6.28318530717959E+0000}%
%
\special{pn 4}%
\special{sh 1}%
\special{ar 2028 1778 16 16 0  6.28318530717959E+0000}%
\special{sh 1}%
\special{ar 2287 1470 16 16 0  6.28318530717959E+0000}%
%
\special{pn 4}%
\special{sh 1}%
\special{ar 1916 1855 16 16 0  6.28318530717959E+0000}%
\special{sh 1}%
\special{ar 1916 1855 16 16 0  6.28318530717959E+0000}%
%
\special{pn 4}%
\special{sh 1}%
\special{ar 1090 1785 16 16 0  6.28318530717959E+0000}%
\special{sh 1}%
\special{ar 1090 1785 16 16 0  6.28318530717959E+0000}%
%
\special{pn 4}%
\special{sh 1}%
\special{ar 1244 1827 16 16 0  6.28318530717959E+0000}%
\special{sh 1}%
\special{ar 1244 1827 16 16 0  6.28318530717959E+0000}%
%
\special{pn 4}%
\special{sh 1}%
\special{ar 3428 770 16 16 0  6.28318530717959E+0000}%
\special{sh 1}%
\special{ar 3435 770 16 16 0  6.28318530717959E+0000}%
%
\special{pn 4}%
\special{sh 1}%
\special{ar 3540 1162 16 16 0  6.28318530717959E+0000}%
\special{sh 1}%
\special{ar 3540 1162 16 16 0  6.28318530717959E+0000}%
%
\special{pn 4}%
\special{sh 1}%
\special{ar 3974 1358 16 16 0  6.28318530717959E+0000}%
\special{sh 1}%
\special{ar 3974 1358 16 16 0  6.28318530717959E+0000}%
\put(14.1000,-22.0000){\makebox(0,0)[lb]{$\xi_1$}}%
\put(35.4000,-21.5000){\makebox(0,0)[lb]{$\xi_2$}}%
%
\special{pn 4}%
\special{sh 1}%
\special{ar 3060 1350 16 16 0  6.28318530717959E+0000}%
\special{sh 1}%
\special{ar 3110 1650 16 16 0  6.28318530717959E+0000}%
\special{sh 1}%
\special{ar 3110 1650 16 16 0  6.28318530717959E+0000}%
%
\special{pn 4}%
\special{sh 1}%
\special{ar 3270 1780 16 16 0  6.28318530717959E+0000}%
\special{sh 1}%
\special{ar 3370 1700 16 16 0  6.28318530717959E+0000}%
\special{sh 1}%
\special{ar 4040 1680 16 16 0  6.28318530717959E+0000}%
\special{sh 1}%
\special{ar 4110 1800 16 16 0  6.28318530717959E+0000}%
\special{sh 1}%
\special{ar 4000 1820 16 16 0  6.28318530717959E+0000}%
\special{sh 1}%
\special{ar 3390 1840 16 16 0  6.28318530717959E+0000}%
\special{sh 1}%
\special{ar 3600 1340 16 16 0  6.28318530717959E+0000}%
\special{sh 1}%
\special{ar 2790 1460 16 16 0  6.28318530717959E+0000}%
\end{picture}}%

\caption{$\xi_1$ is horocyclic and $\xi_2$ is nonhorocyclic.}
\end{figure}

\begin{lemma}\label{l1}
For any lift $\tilde v$ of $v\in T^1\Sigma$, the geodesic ray
 $v[0,\infty)$ in $\Sigma$
is a quasi-minimizer if and only if $\xi=\tilde v(\infty)$
is a nonhorocyclic limit point.
\end{lemma}

\begin{definition} For any point $\xi\in\partial_\infty\HH$, let
$\tilde v\in T^1_i\HH$ be a tangent vector at $i\in\HH$ such that $\tilde v(\infty)=\xi$.
The {\em Buseman function} $B_\xi:\HH\to\R$ is defined for $z\in\HH$ by
\begin{equation}\label{e1}
B_\xi(z)=\lim_{t\to\infty}(d(z, \pi_1\tilde g^t(\tilde v))-t).
\end{equation}
\end{definition}

Notice that for $k>0$, the set $\{B_\xi<-k\}$ is a horodisk at $\xi$ 
which is $k$-apart from $i$.

{\sc Proof of Lemma \ref{l1}}. One may assume that 
$\tilde v$ in the lemma is a unit tangent vector at $i\in\HH$.
Suppose that the point $\xi=\tilde v(\infty)$ is  
a nonhorocyclic limit point. Then there is $k>0$
such that for any $\gamma\in\Gamma$, $B_\xi(\gamma i)\geq -k$.
Since the limit in (\ref{e1}) is non increasing, this implies
$d(\gamma i,\pi_1\tilde g^t(\tilde v))\geq t-k$ for any 
$\gamma\in\Gamma$ and $t\geq0$. On $\Sigma=\Gamma\setminus\HH$, we get
$d(\pi(v),\pi g^t(v))\geq t-k$ for any $t\geq0$. 
That is, $v[0,\infty)$ is a quasi-minimizer.

The converse can be shown by reversing the argument. \qed
@
\begin{lemma}
For any $\xi\in \partial_\infty\HH$ and for any 
$p\in\AAA(\xi)\subset\AAA$, the following conditions are equivalent.

(1) $\Gamma p$ is dense in $\AAA$.

(2) $0\in \overline{\Gamma p}$.

(3) $\xi$ is a horocyclic limit point.
\end{lemma}

\bd (3) $\Rightarrow$ (2): By (3), for any $p\in \AAA(\xi)$, there is
$\gamma\in\Gamma$ such that
$\gamma^{-1} i\in D(p)$. That is, $i\in D(\gamma p)$, namely
$\abs{\gamma p}<1$. Since $p$ is an arbitray point of $\AAA(\xi)$
and since the $\Gamma$-action on $\AAA$ commutes 
with the scalar multiplicaton, this implies (2).

\medskip
(2) $\Rightarrow$ (1): 
For any $\gamma\in\Gamma\setminus\{e\}$, let
$W^u(\gamma)$ be the ray in $\AAA$ corresponding to the eigenspace of $\gamma$
associated to the eigenvalue whose
absolute value is bigger than 1. In other words,
$$
W^u(\gamma)=\{q\in\AAA\mid\abs{\gamma^{-n}q}\to0,\ \ n\to\infty\}.$$
Assume $p$ satisfies (2).
Then we have 
$\overline{\Gamma p}\cap W^u(\gamma)\neq\emptyset$. See Figure 3.
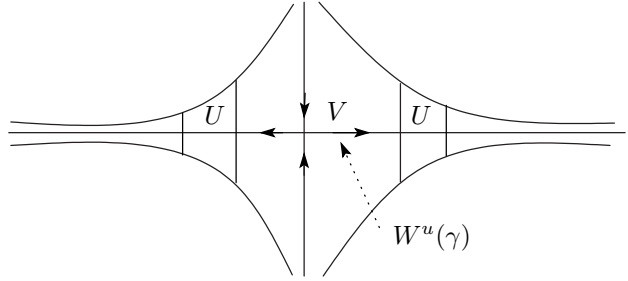
\begin{figure}[h]
{\unitlength 0.1in%
\begin{picture}( 32.2600, 14.3400)(  7.7000,-20.1400)%
%
\special{pn 8}%
\special{pa 770 1263}%
\special{pa 3996 1263}%
\special{fp}%
%
\special{pn 8}%
\special{pa 2316 580}%
\special{pa 2316 2014}%
\special{fp}%
%
\special{pn 8}%
\special{pa 781 1207}%
\special{pa 816 1209}%
\special{pa 852 1212}%
\special{pa 992 1220}%
\special{pa 1027 1221}%
\special{pa 1062 1223}%
\special{pa 1097 1224}%
\special{pa 1131 1224}%
\special{pa 1166 1225}%
\special{pa 1234 1225}%
\special{pa 1267 1224}%
\special{pa 1301 1223}%
\special{pa 1334 1221}%
\special{pa 1366 1219}%
\special{pa 1399 1216}%
\special{pa 1431 1213}%
\special{pa 1463 1209}%
\special{pa 1494 1204}%
\special{pa 1556 1192}%
\special{pa 1616 1178}%
\special{pa 1646 1169}%
\special{pa 1674 1160}%
\special{pa 1703 1149}%
\special{pa 1731 1138}%
\special{pa 1758 1126}%
\special{pa 1785 1113}%
\special{pa 1811 1098}%
\special{pa 1836 1083}%
\special{pa 1861 1067}%
\special{pa 1886 1049}%
\special{pa 1909 1030}%
\special{pa 1933 1011}%
\special{pa 1955 990}%
\special{pa 1999 946}%
\special{pa 2020 923}%
\special{pa 2041 899}%
\special{pa 2062 874}%
\special{pa 2082 849}%
\special{pa 2101 823}%
\special{pa 2121 796}%
\special{pa 2140 770}%
\special{pa 2159 742}%
\special{pa 2178 715}%
\special{pa 2197 687}%
\special{pa 2215 659}%
\special{pa 2234 631}%
\special{pa 2252 603}%
\special{pa 2260 591}%
\special{fp}%
%
\special{pn 8}%
\special{pa 2394 580}%
\special{pa 2416 605}%
\special{pa 2439 630}%
\special{pa 2461 654}%
\special{pa 2483 679}%
\special{pa 2506 704}%
\special{pa 2528 728}%
\special{pa 2574 776}%
\special{pa 2643 845}%
\special{pa 2667 867}%
\special{pa 2690 889}%
\special{pa 2714 910}%
\special{pa 2739 931}%
\special{pa 2763 951}%
\special{pa 2788 971}%
\special{pa 2813 990}%
\special{pa 2838 1008}%
\special{pa 2890 1042}%
\special{pa 2916 1058}%
\special{pa 2943 1073}%
\special{pa 2970 1087}%
\special{pa 2998 1100}%
\special{pa 3026 1112}%
\special{pa 3084 1134}%
\special{pa 3113 1144}%
\special{pa 3142 1153}%
\special{pa 3172 1161}%
\special{pa 3202 1168}%
\special{pa 3233 1175}%
\special{pa 3295 1187}%
\special{pa 3326 1192}%
\special{pa 3390 1200}%
\special{pa 3454 1206}%
\special{pa 3520 1210}%
\special{pa 3552 1212}%
\special{pa 3586 1213}%
\special{pa 3652 1215}%
\special{pa 3820 1215}%
\special{pa 3854 1214}%
\special{pa 3888 1214}%
\special{pa 3922 1213}%
\special{pa 3940 1213}%
\special{fp}%
%
\special{pn 8}%
\special{pa 781 1330}%
\special{pa 817 1328}%
\special{pa 854 1326}%
\special{pa 998 1318}%
\special{pa 1034 1317}%
\special{pa 1070 1315}%
\special{pa 1105 1314}%
\special{pa 1141 1314}%
\special{pa 1176 1313}%
\special{pa 1245 1313}%
\special{pa 1280 1314}%
\special{pa 1314 1316}%
\special{pa 1348 1317}%
\special{pa 1381 1319}%
\special{pa 1414 1322}%
\special{pa 1447 1326}%
\special{pa 1511 1334}%
\special{pa 1573 1346}%
\special{pa 1604 1353}%
\special{pa 1634 1360}%
\special{pa 1692 1378}%
\special{pa 1720 1388}%
\special{pa 1748 1399}%
\special{pa 1775 1411}%
\special{pa 1801 1424}%
\special{pa 1827 1438}%
\special{pa 1852 1453}%
\special{pa 1877 1469}%
\special{pa 1923 1505}%
\special{pa 1945 1524}%
\special{pa 1967 1545}%
\special{pa 1988 1567}%
\special{pa 2008 1590}%
\special{pa 2046 1638}%
\special{pa 2082 1690}%
\special{pa 2099 1717}%
\special{pa 2116 1745}%
\special{pa 2132 1773}%
\special{pa 2149 1802}%
\special{pa 2164 1831}%
\special{pa 2180 1861}%
\special{pa 2195 1890}%
\special{pa 2211 1921}%
\special{pa 2241 1981}%
\special{pa 2254 2008}%
\special{fp}%
%
\special{pn 8}%
\special{pa 2416 2014}%
\special{pa 2435 1986}%
\special{pa 2455 1958}%
\special{pa 2474 1930}%
\special{pa 2494 1903}%
\special{pa 2514 1875}%
\special{pa 2533 1848}%
\special{pa 2553 1821}%
\special{pa 2574 1794}%
\special{pa 2614 1742}%
\special{pa 2677 1667}%
\special{pa 2699 1643}%
\special{pa 2720 1620}%
\special{pa 2743 1598}%
\special{pa 2765 1576}%
\special{pa 2811 1534}%
\special{pa 2859 1496}%
\special{pa 2909 1462}%
\special{pa 2961 1432}%
\special{pa 3015 1406}%
\special{pa 3043 1395}%
\special{pa 3072 1384}%
\special{pa 3100 1375}%
\special{pa 3130 1366}%
\special{pa 3159 1358}%
\special{pa 3189 1351}%
\special{pa 3251 1339}%
\special{pa 3282 1334}%
\special{pa 3313 1330}%
\special{pa 3377 1324}%
\special{pa 3476 1318}%
\special{pa 3509 1318}%
\special{pa 3542 1317}%
\special{pa 3610 1317}%
\special{pa 3746 1321}%
\special{pa 3781 1323}%
\special{pa 3815 1325}%
\special{pa 3850 1326}%
\special{pa 3884 1328}%
\special{pa 3917 1330}%
\special{fp}%
%
\special{pn 8}%
\special{pa 2820 1006}%
\special{pa 2820 1521}%
\special{fp}%
%
\special{pn 8}%
\special{pa 2316 1487}%
\special{pa 2316 1375}%
\special{fp}%
\special{sh 1}%
\special{pa 2316 1375}%
\special{pa 2296 1442}%
\special{pa 2316 1428}%
\special{pa 2336 1442}%
\special{pa 2316 1375}%
\special{fp}%
%
\special{pn 8}%
\special{pa 2204 1263}%
\special{pa 2092 1263}%
\special{fp}%
\special{sh 1}%
\special{pa 2092 1263}%
\special{pa 2159 1283}%
\special{pa 2145 1263}%
\special{pa 2159 1243}%
\special{pa 2092 1263}%
\special{fp}%
%
\special{pn 8}%
\special{pa 2467 1263}%
\special{pa 2652 1263}%
\special{fp}%
\special{sh 1}%
\special{pa 2652 1263}%
\special{pa 2585 1243}%
\special{pa 2599 1263}%
\special{pa 2585 1283}%
\special{pa 2652 1263}%
\special{fp}%
%
\special{pn 8}%
\special{pa 2321 1056}%
\special{pa 2316 1179}%
\special{fp}%
\special{sh 1}%
\special{pa 2316 1179}%
\special{pa 2339 1113}%
\special{pa 2318 1126}%
\special{pa 2299 1112}%
\special{pa 2316 1179}%
\special{fp}%
\put(17.9400,-12.2000){\makebox(0,0)[lb]{$U$}}%
\put(28.7400,-12.2000){\makebox(0,0)[lb]{$U$}}%
\put(24.3400,-11.9600){\makebox(0,0)[lb]{$V$}}%
%
\special{pn 8}%
\special{pa 1680 1160}%
\special{pa 1680 1380}%
\special{fp}%
%
\special{pn 8}%
\special{pa 1960 990}%
\special{pa 1960 1525}%
\special{fp}%
%
\special{pn 8}%
\special{pa 3060 1120}%
\special{pa 3060 1385}%
\special{fp}%
%
\special{pn 8}%
\special{pa 2710 1760}%
\special{pa 2510 1320}%
\special{dt 0.045}%
\special{sh 1}%
\special{pa 2510 1320}%
\special{pa 2519 1389}%
\special{pa 2532 1369}%
\special{pa 2556 1372}%
\special{pa 2510 1320}%
\special{fp}%
\put(27.8000,-18.8000){\makebox(0,0)[lb]{$W^u(\gamma)$}}%
\end{picture}}%

\caption{\small $U$ is a partial fundamental domain for the action of $\gamma$.
If $\overline{\Gamma p}$ intersects $V$, it also intersects
 $U$. The four hyperbolae can be chosen arbitrarily near the axes.}
\end{figure}

Choose $q\in\overline{\Gamma p}\cap W^u(\gamma)$. Let
$\Gamma_n$ be the fundamental group of the subsurface $\Sigma_n$ in
Definition 1.1. The subgroups $\Gamma_n$ are finitely generated with
Cantor limit sets $\Lambda_n$ and form an exausting sequence of subgroups
of $\Gamma$. Moreover $\cup_n\Lambda_n$ is dense in $\partial_\infty\HH$.
 Let
$\AAA_n$ be the inverse image of $\Lambda_n$ by the canonical projection
$\AAA\to\partial_\infty\HH$. We have $\gamma\in\Gamma_n$ and
$q\in\AAA_n$
for any large
$n$.
By the Hedlund theorem \cite{H}, the
$\Gamma_n$ actions on $\AAA_n$ are minimal. In particular,
 $\overline{\Gamma q}\supset\AAA_n$. Since this holds for any large $n$ and since $\overline{\cup_n\AAA_n}=\AAA$, we obtain $\overline{\Gamma q}=\AAA$.
On the other hand, since $q\in\overline{\Gamma p}$,
we have $\overline{\Gamma p}\supset\overline{\Gamma q}$, showing (1).

\medskip
(1) $\Rightarrow$ (3): For any $p\in \AAA(\xi)$, there is $\gamma$ such
      that
$\abs{\gamma^{-1} p}<1$. Then $i\in D(\gamma^{-1} p)$. We thus have
$\Gamma  i\cap
D(p)\neq\emptyset$ for any horodisk $D(p)$ at $\xi$.  \qed

\begin{lemma}
There are horocyclic limit points and nonhorocyclic limit points.
\end{lemma}

\bd Any point in $\partial_\infty\HH$ which is fixed by any
$\gamma\in\Gamma\setminus\{e\}$ is a horocyclic limit point.
To show the second statement, let $D_i$ be the Dirichlet fundamental
domain of $i\in\HH$.
That is, 
$$D_i=\{z\in\HH\mid d(z,i)\leq d(z,\gamma i),\ \forall \gamma\in\Gamma\}.$$
Then any point $\xi$ of $\overline D_i\cap \partial_\infty\HH$ is a
 nonhorocyclic
limit point. In fact, the horodisk $\{B_\xi<0\}$ contains no point of
$\Gamma i$. See Figure 4.
\qed

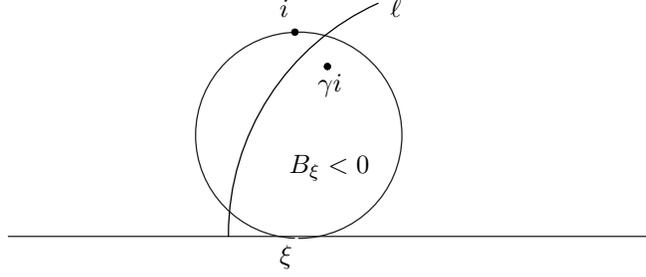
\begin{figure}[h]
{\unitlength 0.1in%
\begin{picture}( 38.0800, 12.8100)(  6.9000,-17.9100)%
%
\special{pn 8}%
\special{pa 690 1590}%
\special{pa 4070 1590}%
\special{fp}%
%
\special{pn 8}%
\special{ar 2210 1060 540 540  1.6124389  1.5707963}%
%
\special{pn 4}%
\special{sh 1}%
\special{ar 2190 520 16 16 0  6.28318530717959E+0000}%
\special{sh 1}%
\special{ar 2190 520 16 16 0  6.28318530717959E+0000}%
%
\special{pn 4}%
\special{sh 1}%
\special{ar 2360 700 16 16 0  6.28318530717959E+0000}%
\special{sh 1}%
\special{ar 2360 700 16 16 0  6.28318530717959E+0000}%
%
\special{pn 8}%
\special{ar 3170 1580 1328 1328  3.1347434  4.2895351}%
\put(21.1000,-4.4000){\makebox(0,0)[lb]{$i$}}%
\put(23.1000,-8.6000){\makebox(0,0)[lb]{$\gamma i$}}%
\put(21.1000,-17.6000){\makebox(0,0)[lb]{$\xi$}}%
\put(26.9000,-4.4000){\makebox(0,0)[lb]{$\ell$}}%
\put(21.6000,-13.0000){\makebox(0,0)[lb]{$B_\xi<0$}}%
\end{picture}}%

\caption{\small If $\gamma i$ is contained in $\{B_\xi<0\}$, then $D_i$ must be
 contained in the region above the perpendicular bisector $\ell$. A contradiction to the definition
of $\xi$.}
\end{figure}

\section{Proof of Theorem \ref{main}}

In this section $\Gamma$ is to be a tight Fuchsian group.

\begin{lemma}
Let $\xi\in \partial_\infty\HH$ be a nonhorocyclic limit point. Then
there is $r>0$ such that for any $p\in \AAA(\xi)$, 
$e^{r/2}p\in\overline{\Gamma p}$.
\end{lemma}

This lemma  implies Theorem \ref{main}. In fact, if $X$ is a minimal set for
the $\Gamma$-action on $\AAA$. Then $X$ must be a proper subset of
$\AAA$ by Lemmata 3.5 and 3.6. Choose $p\in X$ and let
$p\in \AAA(\xi)$. Then $\xi$ is a nonhorocyclic limit point by Lemma 3.5.
The above lemma implies that there is $r>0$ such that
$X\cap e^{r/2}X\neq\emptyset$. Since
$X$ is minimal, this implies $X=e^{r/2}X$, showing that $X$ contains $0$
in its closure. This means that $0\in\overline{\Gamma p}$,
contrary to the fact that $\xi$ is  a nonhorocyclic limit point. 

\bigskip
Lemma 4.1 reduces to the following lemma about the geodesic flow
on $T^1\Sigma$.

\begin{lemma}
Let $\tilde v$ be an arbitrary vector in $T^1\HH$ such that
$\xi=\tilde v(\infty)$ is a nonhorocyclic limit point, and let
$v\in T^1_a\Sigma$ be the projected image of $\tilde v$ ($a\in\Sigma$).
Then there are sequences of vectors $v_n\in T^1_a\Sigma$ and positive
numbers $r_n$ such that $v_n\to v$, $r_n\to r>0$ and
$d(g^{t+r_n}(v_n),g^t(v))\to 0$ as $t\to\infty$.
\end{lemma}

Let us see that Lemma 4.2 implies Lemma 4.1. The last statement shows
that  $g^{r_n}(v_n)$ lies on the strong stable manifold of $v$.
Thus we have $g^{r_n}(v_n)=h^{s_n}(v)$ for some $s_n\in\R$.
We assumed $v_n=g^{-r_n}h^{s_n}(v)\to v$. 
Now the family $\{g^{r_n}\}$ is equiconinuous at $v$, because
$r_n\to r$. Therefore
$d(g^{r_n}(v),h^{s_n}(v))\to0$. That is, $h^{s_n}(v)\to g^r(v)$.
Up on $T^1\HH=PSL(2,\R)$, this means that there are
$\gamma_n\in\Gamma$ such that 
$\gamma_n\tilde h^{s_n}(\tilde v)\to \tilde g^r(\tilde v)$. 
Let $p$ be the projection of $\tilde v$ to $\AAA$.
Then down on $\AAA$, we have $\gamma_np\to e^{r/2}p$, showing Lemma 4.1.

\bigskip 
{\sc Proof of Lemma 4.2.}
It is no loss of generality to assume that $a\in\Sigma_1$, where
$v\in T^1_a\Sigma$. In fact one can take the subsurface $\Sigma_1$
in Definition \ref{def} as large as we want.
By the assumption on $\tilde v$, the
 geodesic ray $v[0,\infty)$ is a quasi-minimizer and thus
proper. Let $t_n$ be the maximum time when $\pi g^{t_n}(v)$
hits $\partial \Sigma_n$. Let $c_n$ be a closed curve on
$\partial\Sigma_n$ starting and ending at $\pi g^{t_n}(v)$.
We choose the direction of $c_n$ in such a way that the tangent vectors
of the curves $v[0,\infty)$ and $c_n$ form an angle $\leq\pi/2$ at
the point $\pi g^{t_n}(v)$. 
See Figure 5.
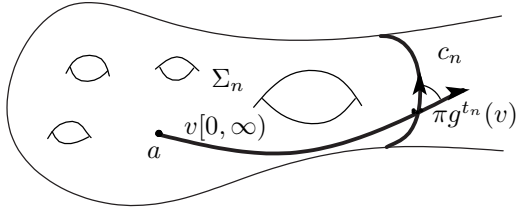
\begin{figure}[h]
{\unitlength 0.1in%
\begin{picture}( 25.9700, 10.6300)(  7.0200,-17.1900)%
%
\special{pn 8}%
\special{pa 3290 726}%
\special{pa 3258 732}%
\special{pa 3227 737}%
\special{pa 3195 743}%
\special{pa 3164 749}%
\special{pa 3132 754}%
\special{pa 3101 760}%
\special{pa 3069 765}%
\special{pa 3038 771}%
\special{pa 3006 776}%
\special{pa 2974 782}%
\special{pa 2943 787}%
\special{pa 2911 792}%
\special{pa 2880 797}%
\special{pa 2848 801}%
\special{pa 2816 806}%
\special{pa 2785 811}%
\special{pa 2721 819}%
\special{pa 2690 823}%
\special{pa 2658 827}%
\special{pa 2626 830}%
\special{pa 2595 834}%
\special{pa 2531 840}%
\special{pa 2499 842}%
\special{pa 2467 845}%
\special{pa 2436 847}%
\special{pa 2404 849}%
\special{pa 2276 853}%
\special{pa 2180 853}%
\special{pa 2116 851}%
\special{pa 2083 849}%
\special{pa 2019 845}%
\special{pa 1987 842}%
\special{pa 1954 839}%
\special{pa 1922 836}%
\special{pa 1890 832}%
\special{pa 1857 827}%
\special{pa 1825 822}%
\special{pa 1792 817}%
\special{pa 1760 811}%
\special{pa 1727 804}%
\special{pa 1695 797}%
\special{pa 1662 790}%
\special{pa 1596 774}%
\special{pa 1564 765}%
\special{pa 1498 745}%
\special{pa 1432 723}%
\special{pa 1399 713}%
\special{pa 1366 702}%
\special{pa 1334 692}%
\special{pa 1301 683}%
\special{pa 1269 674}%
\special{pa 1237 667}%
\special{pa 1205 661}%
\special{pa 1174 658}%
\special{pa 1143 656}%
\special{pa 1113 656}%
\special{pa 1083 659}%
\special{pa 1054 664}%
\special{pa 1025 672}%
\special{pa 997 684}%
\special{pa 969 698}%
\special{pa 943 715}%
\special{pa 917 735}%
\special{pa 892 757}%
\special{pa 869 781}%
\special{pa 846 807}%
\special{pa 825 835}%
\special{pa 805 865}%
\special{pa 787 896}%
\special{pa 770 928}%
\special{pa 755 961}%
\special{pa 741 995}%
\special{pa 730 1030}%
\special{pa 720 1065}%
\special{pa 712 1100}%
\special{pa 707 1136}%
\special{pa 703 1171}%
\special{pa 702 1206}%
\special{pa 703 1241}%
\special{pa 707 1274}%
\special{pa 713 1307}%
\special{pa 721 1339}%
\special{pa 733 1370}%
\special{pa 746 1399}%
\special{pa 762 1428}%
\special{pa 781 1455}%
\special{pa 801 1481}%
\special{pa 823 1505}%
\special{pa 846 1529}%
\special{pa 872 1551}%
\special{pa 899 1571}%
\special{pa 927 1591}%
\special{pa 956 1609}%
\special{pa 986 1626}%
\special{pa 1018 1642}%
\special{pa 1050 1656}%
\special{pa 1082 1668}%
\special{pa 1115 1680}%
\special{pa 1149 1690}%
\special{pa 1182 1698}%
\special{pa 1216 1705}%
\special{pa 1250 1711}%
\special{pa 1283 1715}%
\special{pa 1349 1719}%
\special{pa 1382 1718}%
\special{pa 1414 1717}%
\special{pa 1478 1711}%
\special{pa 1510 1706}%
\special{pa 1541 1700}%
\special{pa 1573 1694}%
\special{pa 1604 1687}%
\special{pa 1635 1679}%
\special{pa 1666 1670}%
\special{pa 1696 1661}%
\special{pa 1758 1641}%
\special{pa 1788 1630}%
\special{pa 1819 1620}%
\special{pa 1849 1609}%
\special{pa 1880 1598}%
\special{pa 1940 1576}%
\special{pa 2002 1554}%
\special{pa 2032 1544}%
\special{pa 2094 1524}%
\special{pa 2125 1515}%
\special{pa 2155 1506}%
\special{pa 2186 1497}%
\special{pa 2218 1489}%
\special{pa 2280 1473}%
\special{pa 2342 1459}%
\special{pa 2374 1453}%
\special{pa 2405 1447}%
\special{pa 2436 1442}%
\special{pa 2468 1436}%
\special{pa 2500 1432}%
\special{pa 2531 1428}%
\special{pa 2563 1424}%
\special{pa 2594 1421}%
\special{pa 2626 1418}%
\special{pa 2722 1412}%
\special{pa 2753 1411}%
\special{pa 2785 1410}%
\special{pa 2881 1410}%
\special{pa 3009 1414}%
\special{pa 3041 1416}%
\special{pa 3074 1417}%
\special{pa 3234 1427}%
\special{pa 3266 1430}%
\special{pa 3299 1432}%
\special{fp}%
%
\special{pn 20}%
\special{pa 2658 831}%
\special{pa 2690 832}%
\special{pa 2721 838}%
\special{pa 2751 851}%
\special{pa 2778 869}%
\special{pa 2800 892}%
\special{pa 2819 917}%
\special{pa 2833 945}%
\special{pa 2845 975}%
\special{pa 2853 1008}%
\special{pa 2859 1041}%
\special{pa 2863 1076}%
\special{pa 2865 1112}%
\special{pa 2864 1148}%
\special{pa 2862 1183}%
\special{pa 2857 1218}%
\special{pa 2849 1252}%
\special{pa 2838 1283}%
\special{pa 2824 1312}%
\special{pa 2807 1338}%
\special{pa 2785 1360}%
\special{pa 2760 1378}%
\special{pa 2731 1393}%
\special{pa 2702 1406}%
\special{pa 2689 1413}%
\special{fp}%
%
\special{pn 20}%
\special{pa 1488 1345}%
\special{pa 1552 1361}%
\special{pa 1584 1368}%
\special{pa 1616 1376}%
\special{pa 1648 1383}%
\special{pa 1680 1391}%
\special{pa 1712 1398}%
\special{pa 1743 1404}%
\special{pa 1775 1411}%
\special{pa 1807 1417}%
\special{pa 1839 1422}%
\special{pa 1870 1427}%
\special{pa 1902 1432}%
\special{pa 1934 1436}%
\special{pa 1965 1439}%
\special{pa 1997 1442}%
\special{pa 2059 1446}%
\special{pa 2122 1446}%
\special{pa 2184 1444}%
\special{pa 2215 1441}%
\special{pa 2246 1437}%
\special{pa 2276 1433}%
\special{pa 2307 1428}%
\special{pa 2338 1422}%
\special{pa 2368 1415}%
\special{pa 2399 1408}%
\special{pa 2459 1392}%
\special{pa 2490 1382}%
\special{pa 2550 1362}%
\special{pa 2610 1340}%
\special{pa 2670 1316}%
\special{pa 2760 1277}%
\special{pa 2850 1235}%
\special{pa 2879 1220}%
\special{pa 2909 1205}%
\special{pa 2939 1191}%
\special{pa 2969 1175}%
\special{pa 2998 1160}%
\special{pa 3058 1130}%
\special{pa 3087 1114}%
\special{pa 3090 1113}%
\special{fp}%
%
\special{pn 4}%
\special{sh 1}%
\special{ar 1497 1341 16 16 0  6.28318530717959E+0000}%
\special{sh 1}%
\special{ar 1497 1341 16 16 0  6.28318530717959E+0000}%
%
\special{pn 20}%
\special{pa 3062 1122}%
\special{pa 3090 1117}%
\special{fp}%
\special{sh 1}%
\special{pa 3090 1117}%
\special{pa 3021 1109}%
\special{pa 3037 1126}%
\special{pa 3028 1148}%
\special{pa 3090 1117}%
\special{fp}%
%
\special{pn 20}%
\special{pa 2862 1068}%
\special{pa 2862 1049}%
\special{fp}%
\special{sh 1}%
\special{pa 2862 1049}%
\special{pa 2842 1116}%
\special{pa 2862 1102}%
\special{pa 2882 1116}%
\special{pa 2862 1049}%
\special{fp}%
%
\special{pn 8}%
\special{pa 1993 1204}%
\special{pa 2015 1178}%
\special{pa 2059 1128}%
\special{pa 2083 1105}%
\special{pa 2107 1084}%
\special{pa 2132 1064}%
\special{pa 2158 1047}%
\special{pa 2186 1033}%
\special{pa 2216 1022}%
\special{pa 2248 1015}%
\special{pa 2281 1011}%
\special{pa 2314 1012}%
\special{pa 2347 1016}%
\special{pa 2379 1024}%
\special{pa 2409 1036}%
\special{pa 2437 1051}%
\special{pa 2462 1071}%
\special{pa 2484 1093}%
\special{pa 2505 1117}%
\special{pa 2543 1171}%
\special{pa 2558 1195}%
\special{fp}%
%
\special{pn 8}%
\special{pa 2039 1159}%
\special{pa 2062 1181}%
\special{pa 2086 1203}%
\special{pa 2111 1223}%
\special{pa 2137 1241}%
\special{pa 2166 1257}%
\special{pa 2197 1270}%
\special{pa 2229 1280}%
\special{pa 2262 1287}%
\special{pa 2296 1290}%
\special{pa 2329 1289}%
\special{pa 2362 1284}%
\special{pa 2392 1275}%
\special{pa 2420 1261}%
\special{pa 2446 1243}%
\special{pa 2470 1222}%
\special{pa 2492 1198}%
\special{pa 2513 1173}%
\special{pa 2525 1159}%
\special{fp}%
%
\special{pn 8}%
\special{pa 1479 1004}%
\special{pa 1507 981}%
\special{pa 1535 960}%
\special{pa 1563 946}%
\special{pa 1591 940}%
\special{pa 1618 945}%
\special{pa 1645 960}%
\special{pa 1672 981}%
\special{pa 1698 1006}%
\special{pa 1706 1013}%
\special{fp}%
%
\special{pn 8}%
\special{pa 1502 995}%
\special{pa 1518 1024}%
\special{pa 1539 1048}%
\special{pa 1568 1061}%
\special{pa 1601 1062}%
\special{pa 1631 1050}%
\special{pa 1655 1029}%
\special{pa 1676 1003}%
\special{pa 1679 999}%
\special{fp}%
%
\special{pn 8}%
\special{pa 1011 1026}%
\special{pa 1027 997}%
\special{pa 1045 970}%
\special{pa 1069 949}%
\special{pa 1100 936}%
\special{pa 1134 932}%
\special{pa 1168 937}%
\special{pa 1196 951}%
\special{pa 1215 975}%
\special{pa 1228 1005}%
\special{pa 1237 1038}%
\special{pa 1238 1040}%
\special{fp}%
%
\special{pn 8}%
\special{pa 1033 995}%
\special{pa 1045 1026}%
\special{pa 1063 1051}%
\special{pa 1091 1065}%
\special{pa 1125 1068}%
\special{pa 1159 1062}%
\special{pa 1189 1048}%
\special{pa 1214 1028}%
\special{pa 1229 1013}%
\special{fp}%
%
\special{pn 8}%
\special{pa 915 1368}%
\special{pa 932 1341}%
\special{pa 952 1316}%
\special{pa 978 1295}%
\special{pa 1011 1282}%
\special{pa 1044 1278}%
\special{pa 1074 1286}%
\special{pa 1098 1306}%
\special{pa 1117 1333}%
\special{pa 1137 1358}%
\special{pa 1142 1363}%
\special{fp}%
%
\special{pn 8}%
\special{pa 938 1331}%
\special{pa 954 1361}%
\special{pa 975 1384}%
\special{pa 1004 1394}%
\special{pa 1037 1394}%
\special{pa 1070 1386}%
\special{pa 1100 1374}%
\special{pa 1128 1360}%
\special{pa 1129 1359}%
\special{fp}%
\put(29.5500,-9.7400){\makebox(0,0)[lb]{$c_n$}}%
\put(14.3600,-14.7100){\makebox(0,0)[lb]{$a$}}%
\put(29.2700,-13.1700){\makebox(0,0)[lb]{$\pi g^{t_n}(v)$}}%
\put(17.7200,-11.2100){\makebox(0,0)[lb]{$\Sigma_n$}}%
%
\special{pn 8}%
\special{ar 2850 1220 134 134  4.7123890  5.8195377}%
%
\special{pn 4}%
\special{sh 1}%
\special{ar 2850 1230 8 8 0  6.28318530717959E+0000}%
\special{sh 1}%
\special{ar 2870 1220 8 8 0  6.28318530717959E+0000}%
\special{sh 1}%
\special{ar 2830 1230 8 8 0  6.28318530717959E+0000}%
\special{sh 1}%
\special{ar 2830 1220 8 8 0  6.28318530717959E+0000}%
\put(16.3000,-13.8000){\makebox(0,0)[lb]{$v[0,\infty)$}}%
\end{picture}}%

\caption{The curves $v[0,\infty)$ and $c_n$.}
\end{figure}

 By Assumption 1.1, there are $0<c<C$ such
that $c\leq\abs{c_n}\leq C$ for any $n$. (If the boundary curve of
$\partial \Sigma_n$ is too short, we choose $c_n$ as its multiple.)
Form a concatenation $\beta_n^T$ of three curves 
$\pi g^t(v)$ ($0\leq
t\leq t_n$), $c_n$ and $\pi g^t(v)$ ($t_n\leq t\leq T$), where $T$ is some
big number. 
Let $\alpha_n^T$ be the geodesic joining $a$ and $\pi g^T(v)$ in
the homotopy class of $\beta_n^T$. If $T\to\infty$, this curve
converges to a geodesic ray $\pi g^t(v_n)$ for some $v_n\in T^1_a\Sigma$.
Moreover the two geodesic rays
$v[0,\infty)$ and $v_n[0,\infty)$  are asymptotic. See Figure 6.
We have $v_n\to v$ by virtue of the bound $C$ in Definition 1.1. See Figure 7.
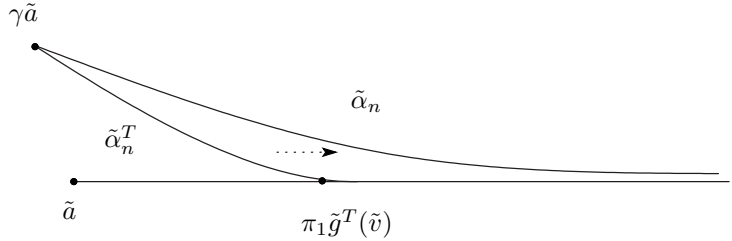
\begin{figure}
{\unitlength 0.1in%
\begin{picture}( 37.7000, 11.1000)(  8.4000,-17.7000)%
%
\special{pn 8}%
\special{pa 1180 1604}%
\special{pa 4610 1604}%
\special{fp}%
%
\special{pn 8}%
\special{pa 970 890}%
\special{pa 1051 941}%
\special{pa 1078 959}%
\special{pa 1213 1044}%
\special{pa 1241 1060}%
\special{pa 1322 1111}%
\special{pa 1350 1127}%
\special{pa 1377 1143}%
\special{pa 1405 1160}%
\special{pa 1432 1176}%
\special{pa 1488 1208}%
\special{pa 1516 1223}%
\special{pa 1572 1255}%
\special{pa 1656 1300}%
\special{pa 1685 1315}%
\special{pa 1713 1329}%
\special{pa 1742 1344}%
\special{pa 1800 1372}%
\special{pa 1829 1385}%
\special{pa 1858 1399}%
\special{pa 1888 1412}%
\special{pa 1917 1425}%
\special{pa 1947 1438}%
\special{pa 2067 1486}%
\special{pa 2098 1497}%
\special{pa 2128 1508}%
\special{pa 2190 1528}%
\special{pa 2252 1546}%
\special{pa 2314 1562}%
\special{pa 2346 1569}%
\special{pa 2377 1576}%
\special{pa 2409 1582}%
\special{pa 2440 1587}%
\special{pa 2472 1592}%
\special{pa 2504 1596}%
\special{pa 2535 1599}%
\special{pa 2567 1602}%
\special{pa 2631 1604}%
\special{pa 2664 1604}%
\special{fp}%
%
\special{pn 8}%
\special{pa 991 897}%
\special{pa 1021 908}%
\special{pa 1051 920}%
\special{pa 1082 931}%
\special{pa 1112 942}%
\special{pa 1142 954}%
\special{pa 1172 965}%
\special{pa 1203 976}%
\special{pa 1233 987}%
\special{pa 1263 999}%
\special{pa 1294 1010}%
\special{pa 1384 1043}%
\special{pa 1415 1054}%
\special{pa 1475 1076}%
\special{pa 1506 1087}%
\special{pa 1536 1098}%
\special{pa 1566 1108}%
\special{pa 1597 1119}%
\special{pa 1627 1130}%
\special{pa 1657 1140}%
\special{pa 1688 1151}%
\special{pa 1718 1161}%
\special{pa 1749 1172}%
\special{pa 1809 1192}%
\special{pa 1840 1202}%
\special{pa 1870 1212}%
\special{pa 1901 1222}%
\special{pa 1931 1232}%
\special{pa 1962 1242}%
\special{pa 1993 1251}%
\special{pa 2023 1261}%
\special{pa 2054 1270}%
\special{pa 2084 1279}%
\special{pa 2146 1297}%
\special{pa 2176 1306}%
\special{pa 2238 1324}%
\special{pa 2269 1332}%
\special{pa 2299 1341}%
\special{pa 2423 1373}%
\special{pa 2454 1380}%
\special{pa 2485 1388}%
\special{pa 2515 1395}%
\special{pa 2546 1402}%
\special{pa 2578 1409}%
\special{pa 2640 1423}%
\special{pa 2764 1447}%
\special{pa 2796 1453}%
\special{pa 2827 1459}%
\special{pa 2858 1464}%
\special{pa 2890 1469}%
\special{pa 2952 1479}%
\special{pa 2984 1483}%
\special{pa 3015 1488}%
\special{pa 3047 1492}%
\special{pa 3078 1496}%
\special{pa 3142 1504}%
\special{pa 3173 1507}%
\special{pa 3205 1511}%
\special{pa 3237 1514}%
\special{pa 3268 1517}%
\special{pa 3364 1526}%
\special{pa 3396 1528}%
\special{pa 3427 1531}%
\special{pa 3619 1543}%
\special{pa 3651 1544}%
\special{pa 3683 1546}%
\special{pa 3715 1547}%
\special{pa 3747 1549}%
\special{pa 3875 1553}%
\special{pa 3908 1554}%
\special{pa 3972 1556}%
\special{pa 4004 1556}%
\special{pa 4068 1558}%
\special{pa 4101 1558}%
\special{pa 4133 1559}%
\special{pa 4197 1559}%
\special{pa 4229 1560}%
\special{pa 4294 1560}%
\special{pa 4326 1561}%
\special{pa 4455 1561}%
\special{pa 4487 1562}%
\special{pa 4554 1562}%
\special{fp}%
%
\special{pn 8}%
\special{pa 2244 1450}%
\special{pa 2552 1450}%
\special{dt 0.045}%
\special{sh 1}%
\special{pa 2552 1450}%
\special{pa 2485 1430}%
\special{pa 2499 1450}%
\special{pa 2485 1470}%
\special{pa 2552 1450}%
\special{fp}%
%
\special{pn 4}%
\special{sh 1}%
\special{ar 977 897 16 16 0  6.28318530717959E+0000}%
\special{sh 1}%
\special{ar 977 897 16 16 0  6.28318530717959E+0000}%
%
\special{pn 4}%
\special{sh 1}%
\special{ar 1180 1604 16 16 0  6.28318530717959E+0000}%
\special{sh 1}%
\special{ar 1180 1604 16 16 0  6.28318530717959E+0000}%
%
\special{pn 4}%
\special{sh 1}%
\special{ar 2480 1600 16 16 0  6.28318530717959E+0000}%
\special{sh 1}%
\special{ar 2480 1600 16 16 0  6.28318530717959E+0000}%
\put(11.2000,-18.0000){\makebox(0,0)[lb]{$\tilde a$}}%
\put(8.4000,-7.9000){\makebox(0,0)[lb]{$\gamma\tilde a$}}%
\put(26.3000,-12.4000){\makebox(0,0)[lb]{$\tilde\alpha_n$}}%
\put(23.7000,-19.0000){\makebox(0,0)[lb]{$\pi_1\tilde g^T(\tilde v)$}}%
\put(13.4000,-14.6000){\makebox(0,0)[lb]{$\tilde \alpha^T_n$}}%
\end{picture}}%

\caption{\small A lift $\tilde v[0,\infty)$ of the curve $v[0,\infty)$
 and a lift $\tilde\alpha_n^T$ of $\alpha_n^T$ to $\HH$. For any $T$, the curves
 $\tilde\alpha_n^T$ start at the same point, say $\gamma\tilde a$, and
converges to $\tilde\alpha_n$. The  projection of $\tilde\alpha_n$ to
$\Sigma$ is the curve $v_n[0,\infty)$.}
\end{figure}
\begin{figure}
{\unitlength 0.1in%
\begin{picture}( 30.4800, 11.3800)(  6.3000,-17.5300)%
%
\special{pn 8}%
\special{pa 630 1374}%
\special{pa 3678 1374}%
\special{fp}%
%
\special{pn 20}%
\special{ar 618 1374 3049 3049  6.0315537  0.1245997}%
%
\special{pn 8}%
\special{pa 1614 1374}%
\special{pa 1866 1032}%
\special{fp}%
%
\special{pn 8}%
\special{pa 1854 1038}%
\special{pa 3618 822}%
\special{fp}%
%
\special{pn 8}%
\special{pa 2322 1374}%
\special{pa 2430 1134}%
\special{fp}%
%
\special{pn 8}%
\special{pa 2430 1146}%
\special{pa 3648 1038}%
\special{fp}%
%
\special{pn 8}%
\special{pa 3054 1374}%
\special{pa 3090 1254}%
\special{fp}%
%
\special{pn 8}%
\special{pa 3096 1260}%
\special{pa 3672 1206}%
\special{fp}%
%
\special{pn 8}%
\special{pa 642 1380}%
\special{pa 3618 822}%
\special{dt 0.045}%
%
\special{pn 8}%
\special{pa 642 1374}%
\special{pa 3642 1038}%
\special{dt 0.045}%
%
\special{pn 8}%
\special{pa 642 1368}%
\special{pa 3660 1206}%
\special{dt 0.045}%
%
\special{pn 4}%
\special{sh 1}%
\special{ar 642 1362 16 16 0  6.28318530717959E+0000}%
\special{sh 1}%
\special{ar 642 1362 16 16 0  6.28318530717959E+0000}%
\put(6.3600,-15.3000){\makebox(0,0)[lb]{$\tilde a$}}%
\end{picture}}%

\caption{\small The dotted lines are the lifts of curves $v_n[0,\infty)$ which
starts at the same point $\tilde a$. This shows
 that $v_n\to v$.}
\end{figure}
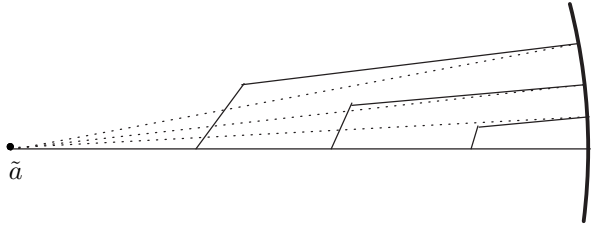

Since 
$v[0,\infty)$ and $v_n[0,\infty)$  are asymptotic, there is $r_n\in\R$ such that
$$d(\pi g^t(v), \pi g^{t+r_n}(v_n))\to 0\ \mbox{ as }\ t\to\infty.$$
 The directions are also 
asymptotic, and therefore 
$$d(g^t(v),g^{t+r_n}(v_n))\to 0\ \mbox{ in }\ T^1\Sigma.$$
Finally we have $r_n\in(b, C]$, where
$b>0$ is defined as follows: for any point $z$ in a horocycle $H$, let $[z,w]$ be the geodesic
segment of length $c/2$ tangent to $H$ at $z$ and pointing outwards. 
Define $b$ by
$b=B_\xi(w)-B_\xi(z)$. For details see Figure 8. This shows Lemma
4.2. \qed

\medskip
The proof of Theorem \ref{main} is now complete. 
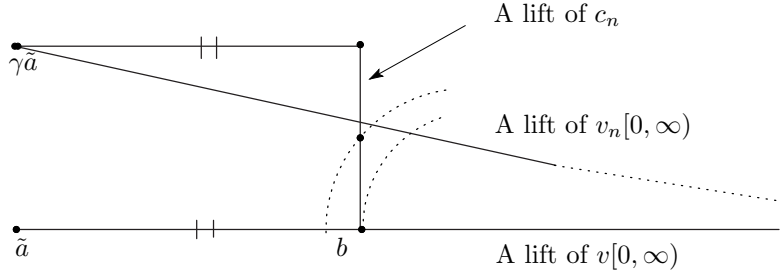
\begin{figure}[h]
{\unitlength 0.1in%
\begin{picture}( 40.3200, 19.2800)(  7.7000,-20.7800)%
%
\special{pn 8}%
\special{pa 794 788}%
\special{pa 2610 788}%
\special{fp}%
%
\special{pn 8}%
\special{pa 802 1748}%
\special{pa 4802 1748}%
\special{fp}%
%
\special{pn 8}%
\special{pa 2610 788}%
\special{pa 2610 1748}%
\special{fp}%
%
\special{pn 8}%
\special{pa 802 788}%
\special{pa 3634 1412}%
\special{fp}%
%
\special{pn 8}%
\special{pa 1778 740}%
\special{pa 1778 852}%
\special{fp}%
%
\special{pn 8}%
\special{pa 1858 740}%
\special{pa 1858 852}%
\special{fp}%
%
\special{pn 8}%
\special{pa 1754 1700}%
\special{pa 1754 1796}%
\special{fp}%
%
\special{pn 8}%
\special{pa 1842 1700}%
\special{pa 1842 1788}%
\special{fp}%
%
\special{pn 4}%
\special{sh 1}%
\special{ar 802 788 16 16 0  6.28318530717959E+0000}%
\special{sh 1}%
\special{ar 818 788 16 16 0  6.28318530717959E+0000}%
%
\special{pn 4}%
\special{sh 1}%
\special{ar 810 1748 16 16 0  6.28318530717959E+0000}%
\special{sh 1}%
\special{ar 810 1748 16 16 0  6.28318530717959E+0000}%
%
\special{pn 4}%
\special{sh 1}%
\special{ar 2610 1268 16 16 0  6.28318530717959E+0000}%
\special{sh 1}%
\special{ar 2610 1268 16 16 0  6.28318530717959E+0000}%
%
\special{pn 4}%
\special{sh 1}%
\special{ar 2618 1748 16 16 0  6.28318530717959E+0000}%
\special{sh 1}%
\special{ar 2618 1748 16 16 0  6.28318530717959E+0000}%
%
\special{pn 8}%
\special{pa 3626 1412}%
\special{pa 4786 1596}%
\special{dt 0.045}%
%
\special{pn 8}%
\special{pn 8}%
\special{pa 2431 1765}%
\special{pa 2431 1757}%
\special{fp}%
\special{pa 2432 1720}%
\special{pa 2432 1713}%
\special{fp}%
\special{pa 2434 1676}%
\special{pa 2435 1668}%
\special{fp}%
\special{pa 2440 1633}%
\special{pa 2441 1625}%
\special{fp}%
\special{pa 2448 1590}%
\special{pa 2450 1582}%
\special{fp}%
\special{pa 2459 1548}%
\special{pa 2461 1540}%
\special{fp}%
\special{pa 2472 1506}%
\special{pa 2474 1498}%
\special{fp}%
\special{pa 2488 1465}%
\special{pa 2490 1457}%
\special{fp}%
\special{pa 2506 1424}%
\special{pa 2510 1417}%
\special{fp}%
\special{pa 2526 1384}%
\special{pa 2531 1378}%
\special{fp}%
\special{pa 2549 1346}%
\special{pa 2554 1340}%
\special{fp}%
\special{pa 2575 1310}%
\special{pa 2579 1304}%
\special{fp}%
\special{pa 2602 1276}%
\special{pa 2606 1270}%
\special{fp}%
\special{pa 2630 1243}%
\special{pa 2636 1237}%
\special{fp}%
\special{pa 2662 1212}%
\special{pa 2667 1206}%
\special{fp}%
\special{pa 2694 1183}%
\special{pa 2700 1178}%
\special{fp}%
\special{pa 2729 1155}%
\special{pa 2735 1151}%
\special{fp}%
\special{pa 2764 1130}%
\special{pa 2771 1126}%
\special{fp}%
\special{pa 2802 1107}%
\special{pa 2809 1103}%
\special{fp}%
\special{pa 2841 1086}%
\special{pa 2849 1083}%
\special{fp}%
\special{pa 2882 1068}%
\special{pa 2889 1064}%
\special{fp}%
\special{pa 2923 1052}%
\special{pa 2930 1049}%
\special{fp}%
\special{pa 2965 1038}%
\special{pa 2972 1036}%
\special{fp}%
\special{pa 3007 1027}%
\special{pa 3015 1026}%
\special{fp}%
\special{pa 3050 1019}%
\special{pa 3058 1018}%
\special{fp}%
%
\special{pn 8}%
\special{pn 8}%
\special{pa 2625 1748}%
\special{pa 2625 1740}%
\special{fp}%
\special{pa 2626 1703}%
\special{pa 2626 1695}%
\special{fp}%
\special{pa 2630 1659}%
\special{pa 2630 1651}%
\special{fp}%
\special{pa 2636 1615}%
\special{pa 2638 1607}%
\special{fp}%
\special{pa 2647 1572}%
\special{pa 2649 1565}%
\special{fp}%
\special{pa 2660 1530}%
\special{pa 2662 1523}%
\special{fp}%
\special{pa 2676 1489}%
\special{pa 2680 1482}%
\special{fp}%
\special{pa 2695 1449}%
\special{pa 2700 1443}%
\special{fp}%
\special{pa 2718 1411}%
\special{pa 2723 1404}%
\special{fp}%
\special{pa 2744 1374}%
\special{pa 2748 1367}%
\special{fp}%
\special{pa 2771 1339}%
\special{pa 2777 1333}%
\special{fp}%
\special{pa 2801 1307}%
\special{pa 2807 1302}%
\special{fp}%
\special{pa 2834 1277}%
\special{pa 2840 1272}%
\special{fp}%
\special{pa 2869 1249}%
\special{pa 2875 1244}%
\special{fp}%
\special{pa 2905 1224}%
\special{pa 2912 1219}%
\special{fp}%
\special{pa 2944 1202}%
\special{pa 2950 1198}%
\special{fp}%
\special{pa 2983 1182}%
\special{pa 2990 1179}%
\special{fp}%
\special{pa 3024 1166}%
\special{pa 3031 1163}%
\special{fp}%
\put(33.1400,-12.7600){\makebox(0,0)[lb]{A lift of $v_n[0,\infty)$}}%
\put(24.8200,-18.8400){\makebox(0,0)[lb]{$b$}}%
\put(7.7000,-9.4000){\makebox(0,0)[lb]{$\gamma\tilde a$}}%
\put(8.0200,-18.8400){\makebox(0,0)[lb]{$\tilde a$}}%
\put(33.2000,-19.6000){\makebox(0,0)[lb]{A lift of $v[0,\infty)$}}%
\put(33.1000,-6.8000){\makebox(0,0)[lb]{A lift of $c_n$}}%
%
\special{pn 4}%
\special{pa 3160 680}%
\special{pa 2660 960}%
\special{fp}%
\special{sh 1}%
\special{pa 2660 960}%
\special{pa 2728 945}%
\special{pa 2707 934}%
\special{pa 2708 910}%
\special{pa 2660 960}%
\special{fp}%
%
\special{pn 4}%
\special{sh 1}%
\special{ar 2610 780 16 16 0  6.28318530717959E+0000}%
\special{sh 1}%
\special{ar 2610 780 16 16 0  6.28318530717959E+0000}%
\end{picture}}%

\caption{The figure depicts the case where the length of $c_n$ takes the minimal value $c$ and $v[0,\infty)$ intersects
 $c_n$ perpendicularly. The dotted circles are horocycles.
Note that the lift of $v_n[0,\infty)$ intersects
the lift of $c_n$ at a point above the midpoint.
Thus we have
 $B_\xi(\gamma\tilde a)- B_\xi(\tilde a)>b$.}
\end{figure}

\section{Examples and remarks}

Let $\Gamma$ be a tight Fuchsian group.
The bi-invariant Haar measure of $PSL(2,\R)$ induces a measure $m$
on $T^1\Sigma=\Gamma\setminus PSL(2,\R)$ invariant both by the geodesic
and horocycle flows. The spaces $\partial_\infty\HH$ and $\AAA$ are
equipped with the standard Lebesgue measures. In what follows, all the
statements concerning the measures are to be with respect to these measures.
The left $\Gamma$ action on $\AAA$ is Morita equivalent to the horocycle
flow on $T^1\Sigma$ in the measure theoretic sense. Especially the
former is ergodic if and only if the latter is (\cite{Z}, 2.2.3).

Denote by $\Lambda_h$ the set of the horocyclic limit points. For a point
$z\in\HH$, denote the Dirichlet fundamental domain of $z$ by $D_z$, i.e.
$$
D_z=\{w\in\HH\mid d(w,z)\leq d(w,\gamma z),\ \ \forall \gamma\in\Gamma\}.$$
Let $F_z=\overline D_z\cap\partial_\infty\HH$ and $E_z=\Gamma F_z$.
Sullivan \cite{Su} showed that $E_z\cup\Lambda_h$ is a full measure set of
$\partial_\infty \HH$ (for any Fuchsian group).

\medskip
{\sc Proof of Corllary \ref{cor}}. By virtue of the Sullivan theorem
and Lemma 3.5, we only need to show that $F_z$ is a null set.
Assume the contrary. For any 
$\gamma\in
\Gamma\setminus\{e\}$, $\gamma F_z\cap F_z$ is at most two points,
lying on the bisector of $z$ and $\gamma z$.
Let 
$$F'_z=F_z\setminus\bigcup_{\gamma\in\Gamma\setminus\{e\}}\gamma F_z,$$
and let $B$ be the inverse image of $F'_z$ by the canonincal projection
$\AAA\to\partial_\infty\HH$. Then $B$ is a partial measurable
fundamental domain for the $\Gamma$ action on $\AAA$, i.e.\ $B$
is positive measured and $B\cap\gamma B=\emptyset$ for any
$\gamma\in\Gamma\setminus\{e\}$. By the Morita equivalence, the
horocycle flow also has a partial fundamental domain: there is
a positive measured set $A\subset T^1\Sigma$ such that 
$A\cap h^nA=\emptyset$ for any $n\in\Z\setminus\{0\}$.
To see this, consider the inverse image $\pi_2^{-1}(B)$
by the canonical projection $\pi_2:PSL(2,\R)\to\AAA$. Clearly
the $\Z$ action  $\tilde h^n$ restricted to $\pi_2^{-1}(B)$ admits
a fundamental domain $A$. On the other hand, $\pi_2^{-1}(B)$ can be
embedded in $T^1\Sigma=\Gamma\setminus PSL(2,\R)$ since
$B\cap\gamma B=\emptyset$ for any
$\gamma\in\Gamma\setminus\{e\}$.

Let us show that almost all points in $A$ has a proper horocycle orbit.
Choose an arbitrary compact set $K$ of $T^1\Sigma$ and let
$$
a_n=m(A\cap h^{-n}(K))=m(h^n(A)\cap K).
$$
Then we have 
$$
\sum_{n\in\Z}a_n\leq m(K)<\infty.
$$
For any $n_0\in\N$, we have
$$
m(A\cap\bigcup_{\abs{n}\geq n_0}h^{-n}(K))\leq\sum_{\abs{n}\geq
n_0}a_n.$$
Therefore
$$
m(A\cap\bigcap_{n_0\in\N}\bigcup_{\abs{n}\geq n_0}h^{-n}(K))=0.$$
Since $K$ is an arbirtrary compact set, this shows that almost all points in $A$ admits a proper horocyclic
orbit. But a proper orbit is a minimal set. This is against Theorem
\ref{main}, completing the proof of Corollary \ref{cor}.

\medskip
Let us discuss Remark \ref{rem} by examples.

\begin{example}
Let $\Gamma_0$ be a cocompact purely hyperbolic Fuchsian group (a
 surface group), and let $\Gamma$ be a nontrivial normal subgroup
of $\Gamma_0$ of infinite index. Then $\Gamma$ is a tight Fuchsian group. 
\end{example}

When $G=\Gamma_0/\Gamma$ is free abelian, then the horocycle flow
on $\Gamma\setminus PSL(2,\R)$ is known to be ergodic \cite{BL}, \cite{LeS}.
On the other hand, if $G$ is nonamenable, there is a nonconstant bounded
harmonic function on the surface $\Sigma=\Gamma\setminus\HH$ \cite{LyS}.
That is, there is a nonconstant bounded $\Gamma$ invariant measurable
function on $\partial_\infty\HH$, and therefore the $\Gamma$ action on
$\partial_\infty\HH$ is not ergodic. This implies that the $\Gamma$
action on $\AAA$ is not ergodic. By the Morita equivalence, the
horocycle flow on $\Gamma\setminus PSL(2,\R)$ is not ergodic.

\medskip
Let $\FF$ be a surface foliation on a compact manifold. If $\FF$ admits
no transverse invariant measures, then there is a continuous
leafwise Riemannian metric of curvature $-1$ \cite{C}. One may ask the
following question.
\begin{question}
 Are generic leaves of $\FF$ either compact, planar, annular or tight?
\end{question}

See \cite{G} for related topics. This is true for the Hirsch foliation
\cite{Hi} and Lie $G$ foliations. For the latter, see \cite{HMM} for the
idea of the proof.

\end{document}